\documentclass[12pt]{amsart}

\usepackage[utf8]{inputenc}
\usepackage[T2A,T1]{fontenc}
\usepackage[english]{babel}

\usepackage{graphicx}

\usepackage{tikz-cd,amssymb}
\usepackage{xcolor}

\newtheoremstyle{otherthm}
{3pt}
{3pt}
{}
{}
{\rmfamily}
{:}
{.5em}
{}

\newtheorem{thm}{Theorem}[section]

\newtheorem{prop}[thm]{Proposition}

\theoremstyle{definition}
\newtheorem{dfn}[thm]{Definition}

\theoremstyle{remark}

\numberwithin{equation}{section}
\makeatletter
\renewenvironment{proof}[1][\proofname]{%
  \par\pushQED{\qed}\normalfont%
  \topsep6\p@\@plus6\p@\relax
  \trivlist\item[\hskip\labelsep\rmfamily#1\@addpunct{.}]%
  \ignorespaces
}{%
  \popQED\endtrivlist\@endpefalse
}
\makeatother

\begin{document}
\title[Photography principle and invariants of manifolds]{Photography principle, data transmission, and invariants of manifolds}
\author{L.Kauffman}
\address{Louis H Kauffman
, Math Dept, UIC, Chicago, IL\\
loukau@gmail.com,
kauffman@uic.edu}
\author{V.O.Manturov}
\address{Vassily O. Manturov, MIPT\\
vomanturov@yandex.ru}
\author{I.M.Nikonov}
\address{Igor M. Nikonov, MSU\\nim@mail.ru}
\author{S.Kim}
\address{Seongjeong Kim, JLU\\kimseongjeong@jlu.edu.cn}

\begin{abstract}
In the present paper we develop the techniques suggested in \cite{ManturovNikonov}
and the photography principle \cite{ManturovWan} for constructing an invariant of 3-manifolds based on Ptolemy relation. We show that a direct implementation of the techniques leads to a trivial invariant and discuss how this approach can be improved to circumvent the difficulties encountered.
\end{abstract}

 \maketitle

AMS MSC 57M25, 57M27,
 51A20, 05E14, 14N20, 51M15

\section{Introduction}

In \cite{MN}, the authors mentioned some deep and unexpected relations between invariants of braids and invariants of 3-manifolds. The clue is that one of the equations we previously used for constructing invariants of braids, namely, the pentagon equation, is {\em ubiquitous}, i.e., appears everywhere else in mathematics, in particular, in constructing
invariants of 3-manifolds \cite{Felikson, FC, FZ1, FZ2}.

In  \cite{ManturovWan} the second named author (V.O.M) developed a very general method, {\em the photography method} for solving the pentagon equation, and other equations. Roughly speaking, if one has several {\em states} (say, triangulations of a pentagon), some {\em data} for each state (e.g., edge lengths or triangle areas) and a {\em rule (formula)} for translating data from one state to another, then after returning to the initial state, we obtain the same data. 

Saying more standardly, this means that the photography method gives a solution to some equations in terms of ``data'' whatever this data means (in geometrical setting one can have lengths, angles, areas, volumes, etc.); even the amount of data may change, so, it is sometimes hard to say what we mean by an ``equation''.

3-manifolds (as well as manifolds of higher dimensions) can be described as equivalence classes of triangulations modulo Pachner moves \cite{Pachner1, Pachner2}. With a triangulation, one can associate some {\em data} and some {\em data transmission laws}: some set of formal variables corresponding to simplices of a given dimension, and the way these lengths are changing under moves.

Taking this together, we can get to an equation (a system of equations)
needed to get a state-sum associated to triangulations, which will be invariant under, say, (2-3)-Pachner moves. Where to take {\em solutions} to these equations?
\begin{figure}
\centering\includegraphics[width=200pt]{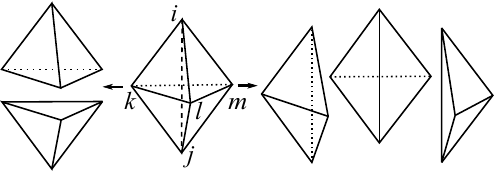}
\caption{(2-3)-Pachner move}
\label{fig:2-3pachner_move}
\end{figure}

First of all, what are the {\em input variables} of the equations? For input variables in 3D-case we may\footnote{Whenever we say ``we may'' we often mean that we may consider variations of formulation, for example using hyperbolic or spherical geometry rather than Euclidean geometry.} take 2-cells, i.e., objects dual to edges of triangulations. The equations will correspond to {\em true vertices} of spines, i.e., objects dual to 3D-simplices. The equations are arranged in such a way that some
{\em state-sum} (or {\em partition function}) is invariant under Pachner moves.

Looking at Pachner moves locally, one sees collections of several tetrahedra
with data (say, lengths) which can be thought of as Euclidean ones\footnote{In fact, the same works perfectly for hyperbolic or spherical ones as well.}. The
Euclidean nature allows one to apply the photography method and guarantees
that some algebraic equations will give rise to an invariant (at the level of
formal variables which originate from Euclidean data).

After that, with each triangulation of a 4-manifold we associate a system
of equations, and show (for free!) that the {\em solution} to this system of
equations (thought of as an algebraic variety) is {\em an invariant} of triangulations under 2-3 or 1-4 moves. Once again, the invariance follows from the {\em Euclidean data transmission law.}

With 1-4 moves one should be more accurate: if we split a tetrahedron
into several small pieces by adding a point, the data for the new triangulation contains more information than the initial data.

In the Section 2 of the paper, we are trying to implement the proposed scheme directly for the Ptolemy relation. The invariant defined formally turns out to be trivial. In Section 3 we discuss possible ways to overcome the problem.


\section{Consistency equation for 3-manifold invariants as a toy model}

Let $T$ be a triangulation of a 3-manifold $M$.
Let $S$ be the 2-frame of the dual cell complex\footnote{To be more accurate, in \cite{ManturovNikonov} we deal with special spines, i.e., graphs dual to the triangular decomposition with only one vertex (which is usually not a triangulation). In this case we need to use Matveev--Piergallini moves instead of Pachner moves.}
corresponding to $T$. Locally, with each 3-simplex of the triangulation, its faces and edges, we associate a true vertex of a complex with four incident edges and six incident 2-cells (a cell is attached to each pair of edges).

Here we shall give the definition of special spine (see~\cite{Matveev}).

\begin{dfn}\label{def:subpolyhedron}
 A \emph{simple polyhedron} is a finite CW-complex $P$ such that each point $x$ in $P$ has a neighbourhood homeomorphic to one of the three configurations (see Figure~\ref{fig:spinePointTypes}):
\begin{enumerate}
   \item a disc;
    \item three discs intersecting along a common arc in their boundaries, where the point $x$ is on the common arc; or
    \item the \emph{butterfly}, which is a configuration built from four arcs that meet at the point $x$, with a face running along each of the six possible pairs of arcs.
\end{enumerate}

\begin{figure}[!ht]
    \centering
    \includegraphics[width=0.25\textwidth]{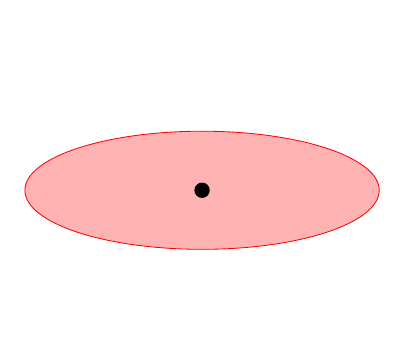}\quad
    \includegraphics[width=0.25\textwidth]{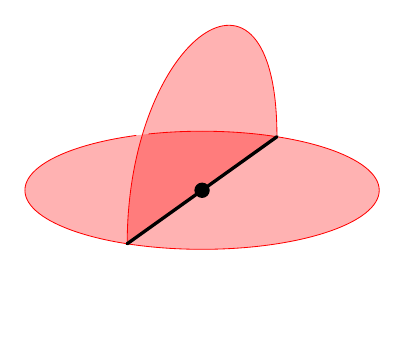}\quad
    \includegraphics[width=0.25\textwidth]{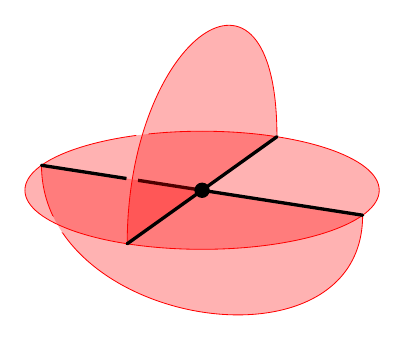}
    \caption{The three types of points in a simple polyhedron $P$.}\label{fig:spinePointTypes}
\end{figure}

\begin{figure}
\centering\includegraphics[width=70pt]{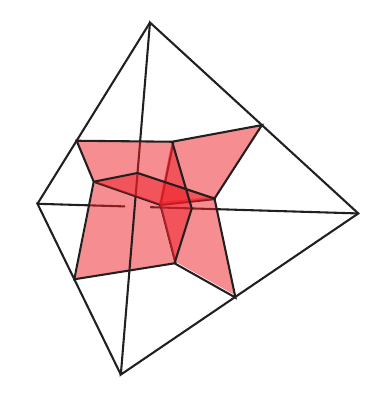}
\caption{Butterfly}
\label{fig:butterfly-tetra}
\end{figure}
Each simple polyhedron is naturally stratified. In this stratification each stratum of dimension $2$ (a $2$-component) is a connected component of the set of nonsingular points. Strata of dimension $1$ consist of open or closed triple lines, and dimension $0$ strata are true vertices.

A simple polyhedron $P$ is called \emph{special} if;
\begin{enumerate}
\item Each 1-stratum of $P$ is an open 1-cell.
\item Each 2-component of $P$ is an open 2-cell.
\end{enumerate}

A \emph{spine} of a compact connected $3$-manifold with boundary $M$ is a subpolyhedron $P\subset Int M$ such that $M\setminus P$ is homeomorphic to $\partial M\times [0,1)$. By a spine of a closed connected $3$-manifold $M$ we mean a spine of $M\setminus Int{} B^3$ where $B^3$ is a $3$-ball in $M$.

\end{dfn}

\subsection{Ptolemy invariant of 3-manifolds}

Let $S$ be a special spine. We associate real-valued variables\footnote{We start with real-valued variables in order to see geometrical reasons why the formula we are going to construct will be invariant; having caught such a formula, we may consider abstract variables.}
to 2-cells of a (special) spine. With each vertex we associate a butterfly; associate this butterfly's weight to it.

Let $S_{0}$ be the set of vertices of the spine $S$, $S_{2}$ the set of 2-cells of $S$, and $\Pi_{v}$ the Ptolemy relation
\begin{equation}\label{eq:Ptolemy_relation}
  \Pi_{v}=ac+bd+xy\in\mathbb Z_2[a,b,c,d,x,y]
\end{equation}
at the vertex $v\in S_0$. Here $a,b,c,d,x,y$ are the weights of the 2-cells incident to $v$, and the the values $a,c$ ($b,d$ and $x,y$) correspond the a pair of opposite cells of the butterfly.
Let $E(S)= \{\Pi_{v}=0, v\in S_{0}\}$ be the set of equations derived from the Ptolemy relations over all true vertices of the special spine.

\begin{figure}
\centering\includegraphics[width=200pt]{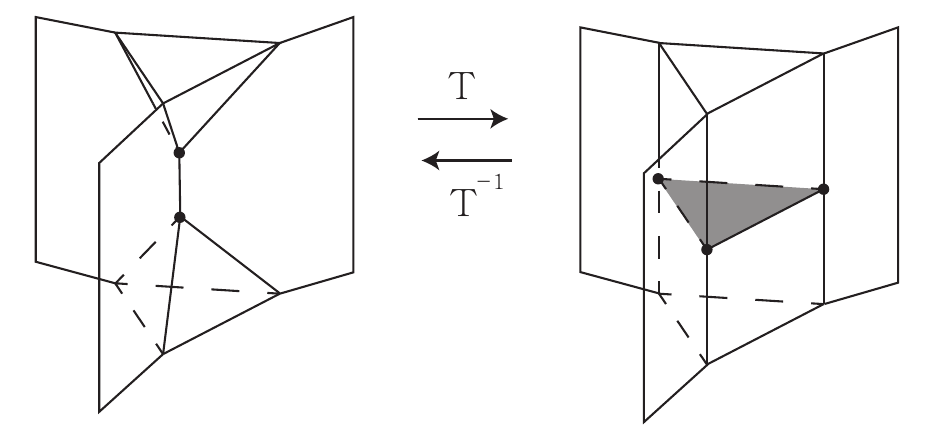}
\caption{The second Matveev--Piergallini move}
\label{fig:Mat-Pill-move-1}
\end{figure}

Let $T,T'$ be two triangulations of $M$ which differ by a 2-3 Pachner move and let $S,S'$ be the corresponding spines. Notice that $S$ and $S'$ differ by a Matveev--Piergallini move (Fig.~\ref{fig:Mat-Pill-move-1}).

Let $E(S), E(S')$ be systems of equations from $S$ and $S'$. It would be pleasant to say that there is a bijection of the solution sets. However, there is an issue concerning division by zero that may occur when one of the weights is zero. 
Let us describe this issue in more details.

Let $K$ be an algebraically closed field of characteristic $2$.

\begin{figure}[h]
\centering\includegraphics[width=0.8\textwidth]{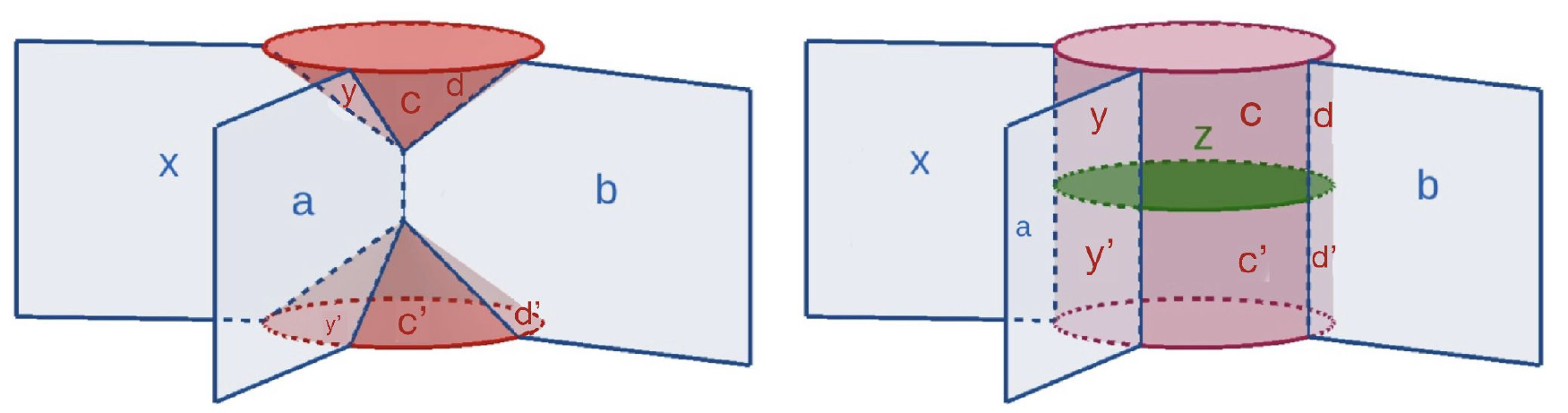}
\caption{The Matveev--Piergallini move}
\label{MP}
\end{figure}

1. Consider a Matveev--Piergallini move  in Fig.\ref{MP}.

Consider the sets of solutions to Ptolemy relations
\[
W=\{(a,b,c,c',d,d',x,y,y')\in K^9\mid\ ad+by+cx=0, ad'+by'+c'x=0\}
\]
and
\begin{multline*}
W'=\{(a,b,c,c',d,d',x,y,y',z)\in K^{10}\mid\\ az+cy'+c'y=0, bz+cd'+c'd=0, xz+dy'+d'y=0\}
\end{multline*}
where $K$ is an algebraically closed field of characteristic $2$.

The projection $p\colon K^{10}\to K^9$, $(a,b,c,c',d,d',x,y,y',z)\mapsto (a,b,c,c',d,d',x,y,y')$ restricts to a map $p|_{W'}\colon W'\to K^9$. Unfortunately, neither $p(W')\subset W$ nor $p(W')\supset W$ nor $p|_{W'}$ is injection. More precisely, we have the following proposition.

\begin{prop}\label{prop:solutions_MP}
1. Let $W_1=\{w\in W\mid a=b=x=0\}$,
\[
W_2=\{w\in W\mid cd'+c'd=0, cy'+c'y=0, dy'+d'y=0\}.
\]
Then
\[
|p^{-1}(w)\cap W'|=\left\{\begin{array}{cl}
1, & w\in W\setminus W_1,\\
0, & w\in W_1\setminus W_2,\\
\infty, & w\in W_1\cap W_2.
\end{array}\right.
\]

2. Let $W_3'=\{w'\in W'\mid z=0\}$. Then $p(W'\setminus W_3')\subset W$ and $p(W'_3)\not\subset W$.
\end{prop}

\begin{proof}
Let $w=(a,b,c,c',d,d',x,y,y')\in W\setminus W_1$. Then either $a\ne 0$ or $b\ne 0$ or $x\ne 0$. Assume $a\ne 0$. Then set $z=\frac{cy'+c'y}a$. Let us check that $w'=(w,z)\in W'$. The equality $az+cy'+c'y=0$ holds by definition of $z$. The equality $bz+cd'+c'd=0$ follows from
\begin{multline*}
a(bz+cd'+c'd)=abz+acd'+ac'd=acd'+ac'd+bcy'+bc'y=\\ c(ad'+by')+c'(ad+by)=cc'x+c'cx=0
\end{multline*}
and $a\ne 0$. Analogously, $xz+dy'+d'y=0$ holds.

Since $z$ is uniquely determined by the equality $az+cy'+c'y=0$, we have $p^{-1}(w)\cap W'=\{w'\}$.

Let $w\in W_1\setminus W_2$. Then, for example, $cy'+c'y\ne 0$. Since $w\in W_1$, $a=0$. Then for any $z\in K$, $az+cy'+c'y=cy'+c'y\ne 0$. Thus, $p^{-1}(w)\cap W'=\emptyset$.

If $w\in W_1\cap W_2$, then the equalities $az+cy'+c'y=0$, $bz+cd'+c'd=0$, $xz+dy'+d'y=0$ hold for any $z\in K$. Hence, $p^{-1}(w)\cap W'\simeq K$.

2. Let $w'\in W'\setminus W'_3$. Then $z\ne 0$ and
\[
a=\frac{cy'+c'y}z,\quad b=\frac{cd'+c'd}z,\quad x=\frac{dy'+d'y}z.
\]
Let us check that $p(w')\in W$.
\[
ad+by+cx = \frac{cdy'+c'dy+cd'y+c'dy+cdy'+cd'y}z=0.
\]
The equality $ad'+by'+c'x=0$ is checked analogously.

Let $w'\in W'_3$. Then $z=0$. Hence, $cd'+c'd=0$, $cy'+c'y=0$, $dy'+d'y=0$, and $a,b,x$ are arbitrary. The conditions on $c,c',d,d',y,y'$ imply that the triples $(c,d,y)$ and $(c',d',y')$ are proportional. If $(c,d,y)=(c',d',y')=(0,0,0)$, then $w=p(w')\in W$ for any $a,b,x$.

Assume that $(c,d,y)\ne(0,0,0)$. Then for any $a,b,x$ such that $ad+by+cz\ne 0$  $p(w')\not\in W$.
\end{proof}

The proposition implies that the maximal domain $\tilde W'\subset W'$ such that $p(\tilde W')\subset W$, $\tilde W'=p^{-1}(p(\tilde W'))\cap W'$ and $p|_{\tilde W'}\colon \tilde W'\to p(\tilde W')$ is a bijection, is equal to
\[
\tilde W'=p^{-1}(W\setminus W_1)\cap W'=W'\setminus W'_1
\]
where $W'_1=\{w'=(a,b,c,c',d,d',x,y,y',z)\in W'\mid a=b=x=0\}$.

\begin{figure}[h]
\centering\includegraphics[width=0.8\textwidth]{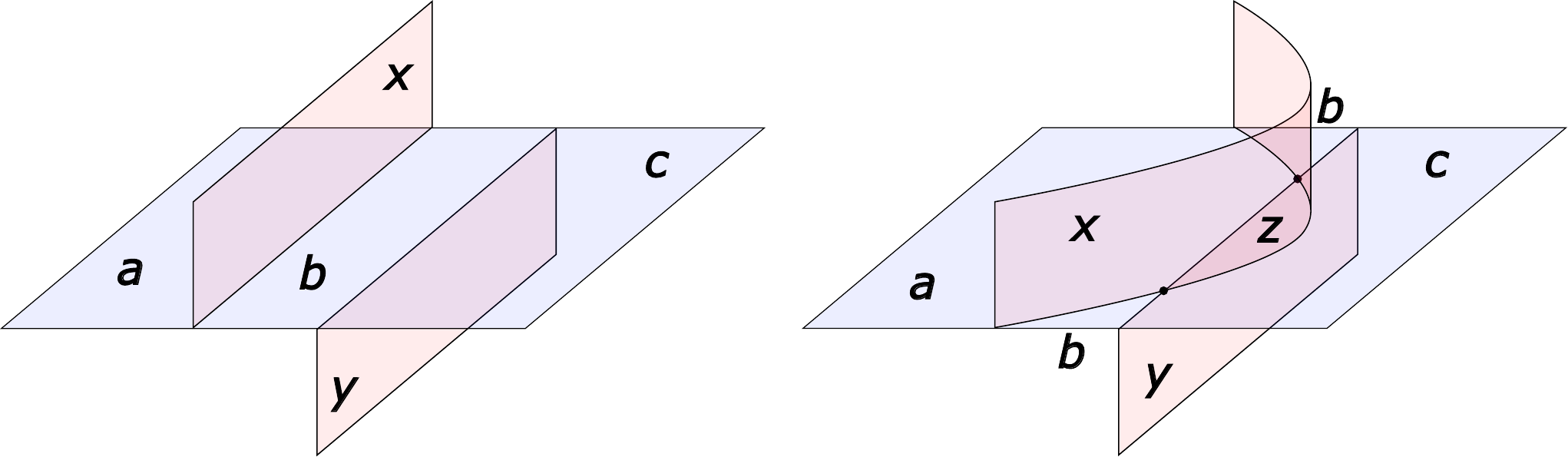}
\caption{The second Matveev--Piergallini move}
\label{MPII}
\end{figure}

2. Consider a second Matveev--Piergallini move (Fig.~\ref{MPII}).

Consider the sets of solutions to Ptolemy relations
$W=\{(a,b,c,x,y)\in K^5\}$
and
\[
W'=\{(a,b,c,x,y,z)\in K^{6}\mid ac+bz+xy=0\}.
\]

 Consider the projection $p\colon K^{6}\to K^5$, $(a,b,c,x,y,z)\mapsto (a,b,c,x,y)$. The following statement is analogous to Proposition~\ref{prop:solutions_MP}.

\begin{prop}\label{prop:solutions_MP2}
1. Let
\[
W_1=\{(a,b,c,x,y)\in K^{5}\mid b= 0\}, \quad W_2=\{(a,b,c,x,y)\in K^{5}\mid ac+xy=0\}.
\]
Then
\[
|p^{-1}(w)\cap W'|=\left\{\begin{array}{cl}
1, & w\in W\setminus W_1,\\
0, & w\in W_1\setminus W_2,\\
\infty, & w\in W_1\cap W_2.
\end{array}\right.
\]
\end{prop}

Thus, the projection $p$ is a bijection between $W\setminus W_1$ and $p^{-1}(W\setminus W_1)\cap W'=W'\setminus W'_1$ where
\[
W'_1=\{(a,b,c,x,y,z)\in W'\mid b= 0\}.
\]

\begin{figure}[h]
\centering\includegraphics[width=0.8\textwidth]{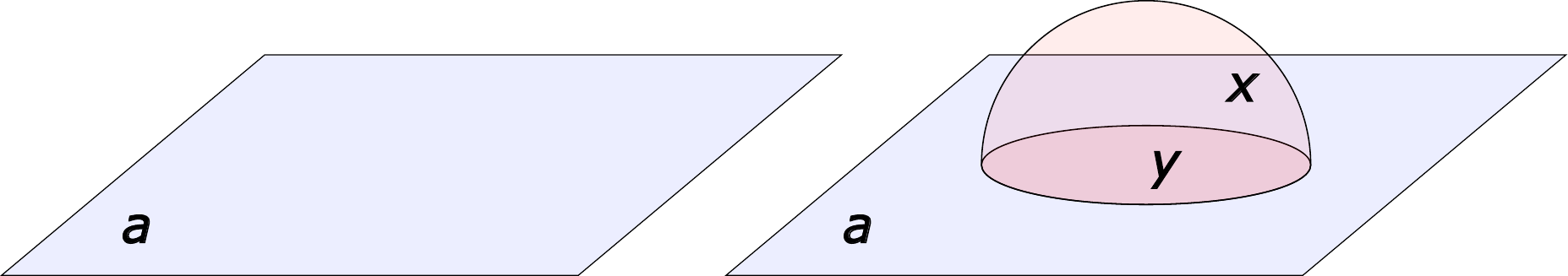}
\caption{The bubble move}
\label{bubble}
\end{figure}

3. Consider the bubble move (Fig.~\ref{bubble}).

The projection $p\colon K^3\to K$, $(a,x,y)\mapsto a$ is a trivial bundle with the fibre $K^2$. Note that $p$ is an open map.

4. Let us construct an invariant of $3$-manifolds, i.e. invariants of spines modulo Matveev--Piergallini moves.

For a spine $S$ with the set of vertices $S_0$ and the set of $2$-cells $S_2$, consider the set of solutions to the Ptolemy relations
\begin{equation}\label{eq:Ptolemy_algebra}
W(S)=\{x=(x_\sigma)_{\sigma\in S_2}\in K^{S_2}\mid \Pi_v(x)=0, v\in S_0\}.
\end{equation}
Denote $W_0(S)=\{(x_\sigma)\in W\mid \forall \sigma\in S_2\ x_\sigma\ne 0\}$. Note that for any Matveev--Piergallini move $f\colon S\to S'$ the induced map $f_*\colon W_0(S)\to W(S')$ is a bijection on the image, and for a bubble map $f_*$ is a trivial bundle with the fibre $K^2$.

Define a sequence of diminishing open subsets $W_k(S)$ in $W(S)$ by the formula
\[
W_{k+1}(S)=W_k(S)\cap\bigcap_{f\colon S\to S'} f_*^{-1}(W_k(S')).
\]
Let $W_\infty(S)=\bigcap_k W_k(S)$.

Let $p\colon K^{n+2}\to K^n$ be a standard projection. We call subsets $W\in K^n$ and $W'\in K^{n+2}$ \emph{stable equivalent} if $p(W')=W'$ and for any $w\in W$ the set $p^{-1}(w)\cap W'$ is dense in $p^{-1}(w)=K^2$.

\begin{thm}
The stable equivalence class of $W_\infty(S)$ is invariant under Matveev--Piergallini moves and bubble moves.
\end{thm}
\begin{proof}
Let $f\colon S\to S'$ is a Matveev--Piergallini move. We show that $f_*\colon W_\infty(S)\to W_\infty(S')$ is a bijection. Indeed, for any $w\in W_\infty(S)$ and any $k$ we have $w\in W_{k+1}(S)\subset f_*^{-1}(W_k(S'))$, hence, $f_*(w)\in W_k(S')$. Then $f_*(w)\in \bigcap_k W_k(S')=W_\infty(S')$ and $f_*(W_\infty(S))\subset W_\infty(S')$. Analogously, $(f_*)^{-1}(W_\infty(S))=(f^{-1})_*(W_\infty(S))\subset W_\infty(S)$. Thus, $f_*$ is a bijection between $W_\infty(S)$ and  $W_\infty(S')$.

Let $f:S\to S'$ is a bubble move. Then $W(S)\subset K^n$ and $W(S')\subset K^{n+2}$, $W(S')=W(S)\times K^2$. The subsets $W_\infty(S)$ and $W_\infty(S')$ are intersections of open subsets $W_k(S)$ in $W(S)$ and $W_k(S')$ in $W(S')$. Since open subsets are dense, the subsets $W_\infty(S)$ and $W_\infty(S')$ are dense in $W(S)$ and $W(S')$. On the other hand, $W_\infty(S')$ projects to $W_\infty(S)$. Hence, $W_\infty(S')$ and $W_\infty(S)$ are stable equivalent.
\end{proof}

\begin{figure}[h]
\centering\includegraphics[width=0.95\textwidth]{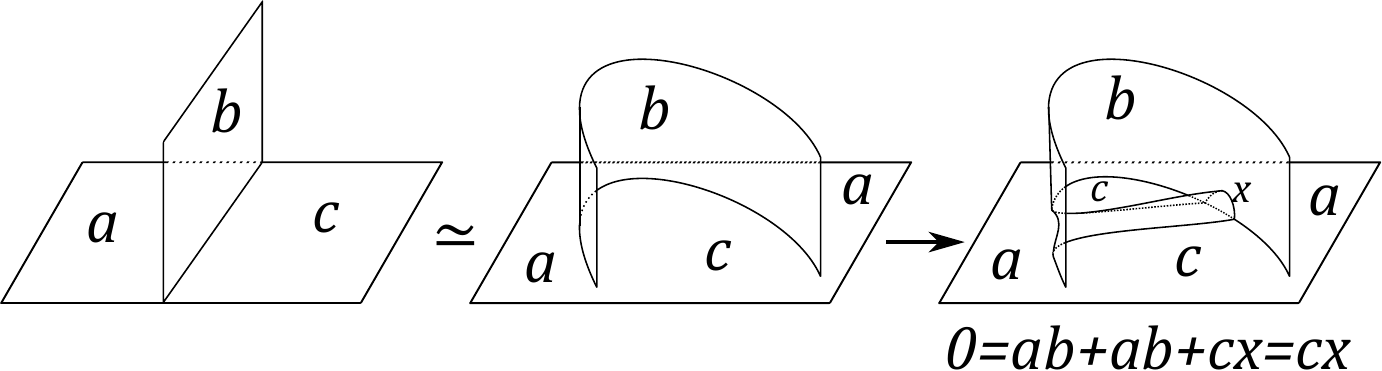}
\caption{Spine transformation}
\label{spine_transformation}
\end{figure}

 Let us show that the set $W_\infty(S)$ is always empty. Indeed, apply a second MP-move $f$ to a spine $S$ as shown in Fig.~\ref{spine_transformation}. Then in the spine after move we have the relation $cx=0$. If $c\ne 0$, then $x=0$. This means that $f(W_0(S))\cap W_0(S')=\emptyset$. Hence, $W_1(S)=f^{-1}(W_0(S'))\cap W_0(S)=\emptyset$. Then $W_\infty(S)=W_1(S)=\emptyset$.

 Thus, we cannot construct a non-trivial invariant while staying in 1-to-1 domain of Matveev--Piergallini moves. If we admit zero 2-faces, then we lose bijective relations between set of solutions to the Ptolemy identity.

\section{Possible difficulties, ways to overcome it, and further direction}

\subsection{Non-uniqueness and non-bijectivity}

 We often relied on the fact: if two triangles $ABCD$, $ABCD'$ (resp., bypiramids $ABCDE$, $ABCDE'$) are both inscribed, then all other bypiramids (say, $ABDEE'$) are inscribed (just because they are both inscribed in the same sphere). This is not true when some of points $A,B,C,D$ coincide.

 In this situation the corresponding system of equations changes but it
 changes in a predictable way. 

Assume we have one system of equations and another system of equations. Well, the solutions between them \emph{are not bijective}. They are \emph{neither injective, nor surjective} but...
They are "almost" bijective provided that \emph{something is non-zero}.
 
But let's write it algebraically. We'll have an algebraic variety on the left and an algebraic
variety on the right. Probably, they will be \emph{birationally equivalent}.
 
I mean, for example, we may have a sphere (well, not a sphere but sphere with some punctured points) on the left and a torus (say, with some singularities) on the right.
 
But we know that \emph{a sphere is not a torus}.

Hence, we should be able to say that "some equivalence class" is \emph{an invariant of our 3-manifold}!

 \subsection{Rigid or mild: why it is invariant and whether it is strong?}

 We present a system of equations for triangulations which are arranged in a way such that their solutions (almost) do not change when Pachner moves are applied.

The reason is that the equations are taken from {\em geometry} (Euclidean, Lobachevskian, etc.) and knowing the data for one triangulation (say, of a pentagon or an octahedron) we can recover the data for the other triangulation.

Let us imagine two inscribed bypiramids sharing the same tetrahedron. Then we can recover the whole polytope with 6 vertices in a rigid way.

This ``rigidity'' is the reason of invariance: all data can be recovered.

However, won't the result be too rigid (roughly speaking, won't we have an arbitrary collection of points corresponding to vertices in the Euclidean space)?

Not quite. First, from geometrical considerations we get equations on lengths, take equations and prove that their spaces of solutions do not change under Pachner moves.

But algebraically, we completely forget about any constraints (triangle inequality etc.) We may work over arbitrary field/ring etc.

Another reason (a much deeper one) is that {\em locally everything is realisable, but globally is not.}

Roughly speaking, we can realise the set of distances corresponding to a triangulation of a simplex in Euclidean (hyperbolic) space, but of course, it is not possible for the whole triangulation is realisable in, say, the Euclidean space, which will lead to some nice equation over the main field; these equation will give rise to an algebra which will be an invariant of our manifold.

\section{Acknowledgements}
We are very grateful to Jim Stasheff, Clifford Henry Taubes and Zheyan Wan for discussion of our current work.

\end{document}